\documentclass[12pt, twoside, leqno]{article}
\usepackage{amsmath,amsthm}
\usepackage{amssymb,latexsym}
\usepackage{enumerate}

\newtheorem{thm}{Theorem}[section]
\newtheorem{prop}[thm]{Proposition}
\newtheorem{cor}[thm]{Corollary}
\newtheorem{lem}[thm]{Lemma}
\newtheorem{prob}[thm]{Problem}

\newtheorem{conj}[thm]{Conjecture}

\theoremstyle{definition}

\newtheorem{rem}[thm]{Remark}

\numberwithin{equation}{section}

\usepackage{mathrsfs}
\usepackage{CJK}
\usepackage{amsmath, amscd}
\usepackage{fancyhdr}
\usepackage{enumerate}

\begin{document}

\baselineskip=17pt

\title{The images of polynomial derivations and $\mathcal{E}$-derivations over field of positive characteristic}
\author{ Fengli Liu \\
MOE-LCSM,\\ School of Mathematics and Statistics,\\
 Hunan Normal University, Changsha 410081, China \\
\emph{E-mail:} liufl4800@163.com
\and
Dan Yan \footnote{ The second author is supported by the NSF of China (Grant No. 12371020), the NSF of Hunan Province (Grant No. 2023JJ30386), the Scientific Research Fund of Hunan Provincial Education Department (Grant No. 21A0056) and the Construct Program of the Key Discipline in Hunan Province.}\\
MOE-LCSM,\\ School of Mathematics and Statistics,\\
 Hunan Normal University, Changsha 410081, China \\
\emph{E-mail:} yan-dan-hi@163.com
}
\date{}

\maketitle

\renewcommand{\thefootnote}{}

\renewcommand{\thefootnote}{\arabic{footnote}}
\setcounter{footnote}{0}

\begin{abstract}
 Let $K$ be a field of characteristic $p$, $\delta$ a nonzero $\mathcal{E}$-derivation and $D=f(x_1)\partial_1$. We first prove that $\operatorname{Im}D$ is not a Mathieu-Zhao space of $K[x_1]$ if and only if $f(x_1)=x_1^rf_1(x_1^p)$ and $r\neq 1$. Then we prove that $\operatorname{Im}\delta$ is a Mathieu-Zhao space of $K[x_1]$ if and only if $\delta$ is not locally nilpotent. Finally, we classify some nilpotent derivations of $K[x_1]$ and give a sufficient and necessary condition for $D(I)$ to be a Mathieu-Zhao space of $K[x_1]$ for any ideal $I$ of $K[x_1]$.
\end{abstract}
{\bf Keywords.} Mathieu-Zhao spaces, $\mathcal{E}$-derivations, Nilpotent derivations.\\
{\bf MSC(2010).} 13N15; 14R10; 13F20; 13G99. \vskip 2.5mm

\section{Introduction}

Throughout this paper, we will write $K$ for a field of characteristic $p$ without specific note, $R$ for an integral domain of characteristic $p$ and $K[x]=K[x_1,x_2,\ldots,x_n]$ ($R[x]=R[x_1,x_2,\ldots,x_n]$) for the
polynomial algebra over $K$ ($R$) in $n$ indeterminates. $\partial_i$ denotes the derivations $\frac{\partial}{\partial x_i}$ for $1\leq i\leq n$, $\operatorname{deg}_{x_i}f$ is the highest degree of $x_i$ in the terms of $f(x_1,\ldots,x_n)$ for $1\leq i\leq n$ and $\operatorname{ldeg}_{x_i}f$ is the least non-negative degree of $x_i$ in the terms of $f(x_1,\ldots,x_n)$ for $1\leq i\leq n$.

An $R$-linear endomorphism $\eta$ of $R[x]$ is said to be locally nilpotent if for each $a\in R[x]$ there exists $m\geq 1$ such that $\eta^m(a)=0$, and locally finite if for each $a\in R[x]$ the $R$-submodule spanned by $\eta^i(a)$ ($i\geq 0$) over $R$ is finitely generated.

A derivation $D$ of $R[x]$ we mean an $R$-linear map $D: R[x]\rightarrow R[x]$ that satisfies $D(ab)=D(a)b+aD(b)$ for all $a, b\in R[x]$ and $D(c)=0$ for any $c\in R$. An $\mathcal{E}$-derivation $\delta$ of $R[x]$ we mean an $R$-linear map $\delta: R[x]\rightarrow R[x]$ such that for all $a, b\in R[x]$ the following equation holds:
$$\delta(ab)=\delta(a)b+a\delta(b)-\delta(a)\delta(b).$$

It is easy to verify that $\delta$ is an $\mathcal{E}$-derivation of $R[x]$, if and only if $\delta=I-\phi$ for some $R$-algebra endomorphism $\phi$ of $R[x]$.

The Mathieu-Zhao space was introduced by Zhao in \cite{1} and \cite{2}, which is a natural generalization of ideals. We give the definition here for the polynomial rings. An $R$-submodule $M$ of $R[x]$ is said to be a Mathieu-Zhao space if for any $a, b\in R[x]$ with $a^m \in M$ for all $m\geq 1$, we have $ba^m\in M$ when $m>>0$. The radical of a Mathieu-Zhao space was first introduced in \cite{2}, denoted by $\mathfrak{r}(M)$, and
$$\mathfrak{r}(M)=\{a\in R[x]| a^m\in M~\operatorname{for}~\operatorname{all}~m>>0\}.$$

There is an equivalent definition about Mathieu-Zhao space which proved in Proposition 2.1 of \cite{2}. We only give the equivalent definition here for the polynomial rings. An $R$-submodule $M$ of $R[x]$ is said to be a Mathieu-Zhao space if for any $a, b\in R[x]$ with $a\in \mathfrak{r}(M)$, we have $ba^m\in M$ when $m>>0$.

In \cite{3}, Wenhua Zhao posed the following two conjectures:
\begin{conj}(LFED)\label{conj1.1}
Let $K$ be a field of characteristic zero and $\mathcal{A}$ a $K$-algebra. Then for every locally finite derivation or $\mathcal{E}$-derivation $\delta$ of $\mathcal{A}$, the image $\operatorname{Im}\delta:=\delta(\mathcal{A})$ of $\delta$ is a Mathieu-Zhao space of $\mathcal{A}$.
\end{conj}

\begin{conj}(LNED)\label{conj1.2}
Let $K$ be a field of characteristic zero and $\mathcal{A}$ a $K$-algebra and $\delta$ a locally nilpotent derivation or $\mathcal{E}$-derivation of $\mathcal{A}$, Then for every ideal $I$ of $\mathcal{A}$, the image $\delta(I)$ of $I$ under $\delta$ is a Mathieu-Zhao space of $\mathcal{A}$.
\end{conj}

There are many positive answers to the above two conjectures.  In \cite{4}, Wenhua Zhao proved that Conjecture \ref{conj1.1} is true for polynomial algebras in one variable and Conjecture \ref{conj1.2} is true for polynomial algebras in one variable for derivations and most $\mathcal{E}$-derivations. Arno van den Essen, David Wright, Wenhua Zhao showed that Conjecture \ref{conj1.1} is true for derivations for polynomial algebras in two variables in \cite{5}. In \cite{6}, Wenhua Zhao proved that Conjecture \ref{conj1.1} is true for Laurent polynomial algebras in one or two variables and Conjecture \ref{conj1.2} is true for all Laurent polynomial algebras. Wenhua Zhao proved the above two conjectures for algebraic algebras in \cite{7}. In \cite{8}, Dayan Liu, Xiaosong Sun showed that Conjecture \ref{conj1.1} is true for linear locally nilpotent derivations in dimension three. They also proved that Conjecture \ref{conj1.1} is true for triangular derivations and homogeneous locally nilpotent derivations in dimension three in \cite{10}. Arno van den Essen, Wenhua Zhao showed that Conjecture \ref{conj1.1} is true for locally integral domains and $K[[x]][x^{-1}]$ in \cite{9}. However, the above two conjectures are not true if we replace a field of characteristic zero by a field of positive characteristic (Example 2.5, Example 2.6 in \cite{3}). Thus, we pose the following problems for the polynomial algebras over a field of positive characteristic.
\begin{prob}\label{prob1.3}
Let $K$ be a field of characteristic $p$ and $\delta$ a derivation or an $\mathcal{E}$-derivation of $K[x]$. Under what conditions on $\delta$ such that $\operatorname{Im}\delta$ is a Mathieu-Zhao space of $K[x]$?
\end{prob}

\begin{prob}\label{prob1.4}
Let $K$ be a field of characteristic $p$, $\delta$ a derivation or an $\mathcal{E}$-derivation of $K[x]$ and $I$ any ideal of $K[x]$. Under what conditions on $\delta$ and $I$ such that $\delta(I)$ is a Mathieu-Zhao space of $K[x]$?
\end{prob}

\begin{prob}\label{prob1.5}
Let $K$ be a field of characteristic $p$ and $\delta$ a derivation or an $\mathcal{E}$-derivation of $K[x]$. Under what conditions on $\delta$ such that $\delta$ is locally nilpotent or locally finite?
\end{prob}

Problem \ref{prob1.3} and Problem \ref{prob1.4} are specific cases of Problem 1.4 $(B)$ and $(C)$ in \cite{2}, respectively. Arrangement is as follows. In section 2, we solve Problem \ref{prob1.3} for $n=1$ and we also prove that the image of a triangular derivation of $K[x]$ is a Mathieu-Zhao space of $K[x]$ if and only if $D=0$. In section 3, we give partial answer to Problem \ref{prob1.5} for $n=1$. Then we solve Problem \ref{prob1.4} for derivations of $K[x_1]$ and give partial answer to Problem \ref{prob1.4} for $\mathcal{E}$-derivations of $K[x_1]$ in section 4.

\section{Problem \ref{prob1.3} for the case $n=1$}

\begin{thm} \label{thm2.1}
Let $D=f(x_1,\ldots,x_n)\partial_{i_0}$ be a derivation of $R[x]$ for some $i_0\in \{1,2,\ldots,n\}$. If $p|\operatorname{deg}_{x_{i_0}}f-1$ or $p|\operatorname{ldeg}_{x_{i_0}}f-1$, then $\mathfrak{r}(\operatorname{Im}D)=0$. In particular, $\operatorname{Im}D$ is a Mathieu-Zhao space of $R[x]$.
\end{thm}
\begin{proof}
Without loss of generality, we can assume that $i_0=1$. Let $\bar{n}=\operatorname{deg}_{x_1}f$ and $\tilde{n}=\operatorname{ldeg}_{x_1}f$. Then we have that
 $$\operatorname{Im}D=R\{x_2^{i_2}\cdots x_n^{i_n}f(x),~x_2^{i_2}\cdots x_n^{i_n}f(x)x_1,\ldots,x_2^{i_2}\cdots x_n^{i_n}f(x)x_1^{p-2},$$
 $$x_2^{i_2}\cdots x_n^{i_n}f(x)x_1^p,\ldots,x_2^{i_2}\cdots x_n^{i_n}f(x)x_1^{kp+t}~\operatorname{for}~\operatorname{all}~0\leq t\leq p-2,~k, ~i_2,\ldots,i_n\in \mathbb{N}\}.$$
 Thus, the polynomials with degree $kp+\bar{n}-1$ of $x_1$ or with least degree $kp+\tilde{n}-1$ of $x_1$ are not in $\operatorname{Im}D$ for all $k\in \mathbb{N}^*$. Suppose that $P(x)$ is a nonzero polynomial in $\mathfrak{r}(\operatorname{Im}D)$. Then we have $P(x)^m\in \operatorname{Im}D$ for all $m\geq N$ for some $N\in \mathbb{N}^*$. Let $g:=P(x)^N=a_i(x_2,\ldots,x_n)x_1^if(x)+\cdots+a_j(x_2,\ldots,x_n)x_1^jf(x)$ for some $a_i,\ldots,a_j\in R[x_2,\ldots,x_n]$, $0\leq i<j$, $a_ia_j\neq 0$.

 $(1)$ If $p|\bar{n}-1$, then $\bar{n}-1=q_1p$ for some $q_1\in \mathbb{N}$. Thus, we have $\operatorname{deg}_{x_1}g^{(q_1+k)p}=(j+\bar{n})(q_1+k)p=[k(j+\bar{n})+q_1(j+\bar{n}-1)]p+\bar{n}-1$. Therefore, we have $g^{(q_1+k)p}\notin \operatorname{Im}D$ for all $k\in \mathbb{N}^*$, which is a contradiction. Thus, we have $P(x)=0$. That is, $\mathfrak{r}(\operatorname{Im}D)=0$.

 $(2)$ If $p|\tilde{n}-1$, then $\tilde{n}-1=q_2p$ for some $q_2\in \mathbb{N}$. Thus, we have $\operatorname{ldeg}_{x_1}g^{(q_2+k)p}=(i+\tilde{n})(q_2+k)p=[k(i+\tilde{n})+q_2(i+\tilde{n}-1)]p+\tilde{n}-1$. Therefore, we have $g^{(q_2+k)p}\notin \operatorname{Im}D$ for all $k\in \mathbb{N}^*$, which is a contradiction. Thus, we have $P(x)=0$. That is, $\mathfrak{r}(\operatorname{Im}D)=0$.
\end{proof}

\begin{cor}\label{cor2.2}
Let $D=f(x_1)\partial_1$ be a derivation of $R[x_1]$. If $p|\operatorname{deg}f(x_1)-1$ or $p|\operatorname{ldeg}f(x_1)-1$, then $\mathfrak{r}(\operatorname{Im}D)=0$. In particular, $\operatorname{Im}D$ is a Mathieu-Zhao space of $R[x_1]$.
\end{cor}
\begin{proof}
 The conclusion follows from Theorem \ref{thm2.1} for $n=1$.
\end{proof}

\begin{prop}\label{prop2.3}
Let $D=x_1^rf_1(x_1^p)\partial_1$ be a nonzero derivation of $R[x_1]$ for $f_1(x_1^p)\in K[x_1^p]$. Then $\operatorname{Im}D$ is a Mathieu-Zhao space of $R[x_1]$ if and only if $r=lp+1$ for some $l\in \mathbb{N}$.
\end{prop}
\begin{proof}
$``\Leftarrow"$ If $r=lp+1$, then $\operatorname{deg}(x_1^rf_1(x_1^p))=r+Lp=(l+L)p+1$, where $Lp=\operatorname{deg}_{x_1}f_1(x_1^p)$. Thus, we have $p|\operatorname{deg}(x_1^rf_1(x_1^p))-1$. It follows from Corollary \ref{cor2.2} that $\mathfrak{r}(\operatorname{Im}D)=0$. Thus, $\operatorname{Im}D$ is a Mathieu-Zhao space of $R[x_1]$.

$``\Rightarrow"$ It suffices to prove that if $r\neq lp+1$ for some $l\in \mathbb{N}$, then $\operatorname{Im}D$ is not a Mathieu-Zhao space of $R[x_1]$. Since $D(x_1^{kp+\tilde{t}})=\tilde{t}x_1^{r+kp+\tilde{t}-1}f_1(x_1^p)$ for $1\leq \tilde{t}\leq p-1$, we have that
$$\operatorname{Im}D=R\{x_1^rf_1(x_1^p),~x_1^{r+1}f_1(x_1^p),\ldots,x_1^{r+p-2}f_1(x_1^p),~x_1^{r+p}f_1(x_1^p),\ldots,x_1^{r+kp+\tilde{t}}f_1(x_1^p)$$
$$ \operatorname{for}~\operatorname{all}~0\leq \tilde{t}\leq p-2,~k\in\mathbb{N}\}.$$
Let $r=lp+t$ for $t=0, 2,\ldots,p-1$, $l\in \mathbb{N}$. Then it's easy to check that $x_1^{kp}f_1(x_1^p)\in \operatorname{Im}D$ for all $k\in \mathbb{N}^*$. In particular, $x_1^pf_1(x_1^p)\in \operatorname{Im}D$. Since
$$(x_1^pf_1(x_1^p))^m=[x_1^{(m-1)p}f_1(x_1^p)^{m-1}]x_1^pf_1(x_1^p)$$
and $x_1^{(m-1)p}f_1(x_1^p)^{m-1}$ is a polynomial of $x_1^p$ for all $m\in \mathbb{N}^*$, we have $(x_1^pf_1(x_1^p))^m\in \operatorname{Im}D$ for all $m\in \mathbb{N}^*$. It's easy to see that $x_1^{r+p-1}(x_1^pf_1(x_1^p))^m\notin \operatorname{Im}D$ for all $m\in \mathbb{N}^*$. Hence, $\operatorname{Im}D$ is not a Mathieu-Zhao space of $R[x_1]$.
\end{proof}

\begin{thm}\label{thm2.4}
Let $D=f(x_1)\partial_1$ be a nonzero derivation of $R[x_1]$. Then we have the following statements:

$(1)$ If $f(x_1)=x_1^rf_1(x_1^p)$ for some $r\in \{0,1,2,\ldots,p-1\}$ and $f_1(x_1^p)\in R[x_1^p]$, then $\operatorname{Im}D$ is a Mathieu-Zhao space of $R[x_1]$ if and only if $r=1$.

$(2)$ If $f(x_1)\neq x_1^rf_1(x_1^p)$ for any $r\in \{0,1,2,\ldots,p-1\}$ and $f_1(x_1^p)\in R[x_1^p]$, then $\mathfrak{r}(\operatorname{Im}D)=0$. In particular, $\operatorname{Im}D$ is a Mathieu-Zhao space of $R[x_1]$.
\end{thm}
\begin{proof}
$(1)$ The conclusion follows from Proposition \ref{prop2.3}.

$(2)$ Since $f(x_1)\in R[x_1]$, we have
$$f(x_1)=a_1x_1^{i_1}f_1(x_1^p)+a_2x_1^{i_2}f_2(x_1^p)+\cdots+a_jx_1^{i_j}f_j(x_1^p)$$
for $a_1,\ldots,a_j\in R^*$, $f_1,\ldots,f_j\in R[x_1^p]$, $f_1\cdots f_j\neq 0$, $j\geq 2$, where $i_1,\ldots,i_j$ are distinct integers between 0 and $p-1$. Thus, we have
$$D(x_1^{kp+t})=t(a_1x_1^{i_1+kp+t-1}f_1(x_1^p)+a_2x_1^{i_2+kp+t-1}f_2(x_1^p)+\cdots+a_jx_1^{i_j+kp+t-1}f_j(x_1^p))$$
for all $1\leq t\leq p-1$, $k\in \mathbb{N}$. Hence, if $\tilde{h}\in \operatorname{Im}D$, then $\tilde{h}$ contains at least two different monomials by module $x_1^p$. We claim that $\mathfrak{r}(\operatorname{Im}D)=0$. Assume otherwise, if there exists a nonzero polynomial $q(x_1)\in \mathfrak{r}(\operatorname{Im}D)$, then we can write
$$q(x_1)=b_1x_1^{j_1}g_1(x_1^p)+b_2x_1^{j_2}g_2(x_1^p)+\cdots+b_ix_1^{j_i}g_i(x_1^p)$$
for $b_1,\ldots,b_i\in R^*$, $g_1,\ldots,g_i\in R[x_1^p]$, $g_1\cdots g_i\neq 0$, where $j_1,\ldots,j_i$ are distinct integers between 0 and $p-1$. Thus, we have
$$q(x_1)^{p^m}=b_1^{p^m}x_1^{j_1p^m}g_1(x_1^p)^{p^m}+\cdots+b_i^{p^m}x_1^{j_ip^m}g_i(x_1^p)^{p^m}=\bar{c}~ \operatorname{mod}~x_1^p$$
for some $\bar{c}\in R$ and for all $m\in \mathbb{N}^*$. Whence we have $q(x_1)\notin \mathfrak{r}(\operatorname{Im}D)$, which is a contradiction. Thus, we have $\mathfrak{r}(\operatorname{Im}D)=0$.
\end{proof}

\begin{thm}\label{thm2.5}
Let $\delta=I-\phi$ be a nonzero $\mathcal{E}$-derivation of $K[x_1]$ and $\phi(x_1)=x_1+c$ for some $c\in K$. Then $\operatorname{Im}\delta$ is not a Mathieu-Zhao space of $K[x_1]$.
\end{thm}
\begin{proof}
Since $\delta$ is nonzero, we have $c\neq 0$. Since $\delta(x_1)=-c\neq 0$, we have $1\in \operatorname{Im}\delta$.
We claim that $x_1^{p-1}\notin \operatorname{Im}\delta$. Suppose that $x_1^{p-1}\in \operatorname{Im}\delta$. Then there exists $u(x_1)\in K[x_1]$ such that
\begin{equation}\label{eq2.1}
x_1^{p-1} = \delta(u(x_1))=u(x_1)-u(x_1+c)
\end{equation}
We have the following equations
\begin{equation}
  \left\{ \begin{aligned}
  \nonumber
 x_1^{p-1} = u(x_1)-u(x_1+c)~~~~~~~~~~~~~~~~~~~~~~~~~~~~~~~~\\
(x_1+c)^{p-1} = u(x_1+c)-u(x_1+2c)~~~~~~~~~~~~~~~~~~\\
\vdots~~~~~~~~~~~~~~~~~~~~~~~~~~~~~~~~~~~~~~~~~\\
(x_1+(p-1)c)^{p-1} = u(x_1+(p-1)c)-u(x_1+pc)
                          \end{aligned} \right.
  \end{equation}
by substituting $x_1$ with $x_1+jc$ for $0\leq j\leq p-1$ in equation \eqref{eq2.1}. Then we have
\begin{equation}\label{eq4.2}
\sum_{j=0}^{p-1}(x_1+jc)^{p-1} = 0
\end{equation}
by adding the above $p$ equations. By Fermat's little Theorem, we have $j^{p-1}=1 \operatorname{mod} p$ for all $1\leq j\leq p-1$. Let $x_1=0$ in equation \eqref{eq4.2}. Then we have that $(p-1)c^{p-1}=0$, which is a contradiction. Hence we have $x_1^{p-1}\notin \operatorname{Im}\delta$. Whence, $\operatorname{Im}\delta$ is not a Mathieu-Zhao space of $K[x_1]$.
\end{proof}

\begin{cor}\label{cor2.6}
 Let $\delta$ be a nonzero $\mathcal{E}$-derivation of $K[x_1]$. Then $\operatorname{Im}\delta$ is a Mathieu-Zhao space of $K[x_1]$ if and only if $\delta$ is not locally nilpotent.
\end{cor}
\begin{proof}
It follows from the discussion after Lemma 3.1 in \cite{4} that $\delta=I-\phi$ can be divided into the following four(exhausting) cases: $\phi(x_1)=c$ for some $c\in K$ or $\phi(x_1)=x_1+c$ for some $c\in K$ or $\phi(x_1)=\tilde{q}x_1$ for some nonzero $\tilde{q}\in K$ or $\operatorname{deg}\phi(x_1)\geq 2$. It's easy to check that $\delta$ is locally nilpotent if and only if $\phi(x_1)=x_1+c$ or $\phi(x_1)=\tilde{q}x_1$ and $\tilde{q}=1$. Since $\delta$ is nonzero, we have that $\delta$ is locally nilpotent if and only if $\phi(x_1)=x_1+c$ and $c\neq 0$. Then the conclusion follows from Lemma 3.2, Lemma 3.3 and Proposition 3.7 in \cite{4} and Theorem \ref{thm2.5}.
\end{proof}

\begin{rem}\label{rem2.7}
Let $\delta$ be as in Corollary \ref{cor2.6}. It's easy to check that $\delta$ is locally finite if and only if $\phi$ is an affine endomorphism.
\end{rem}

\begin{prop}\label{prop2.8}
Let $D=\sum_{q=1}^nf_q(x_{q+1},\ldots,x_n)\partial_q$ be a derivation of $K[x]$, where $f_q(x_{q+1},\ldots,x_n)\in K[x_{q+1},\ldots,x_n]$ for $1\leq q\leq n-1$, $f_n\in K$. Then $\operatorname{Im}D$ is a Mathieu-Zhao space of $K[x]$ if and only if $D=0$.
\end{prop}
\begin{proof}
$``\Leftarrow"$ Clearly.

$``\Rightarrow"$ We assume that $D\neq 0$. Let $g$ be any element in $K[x]$. Then $g$ can be written
$$g=\sum_{\mbox {\tiny $\begin{array}{c}
0\leq i_1\leq p-1\\
\vdots\\
0\leq i_n\leq p-1\\
i_1+\cdots+i_n\geq 1\end{array}$}}x_1^{i_1}x_2^{i_2}\cdots x_n^{i_n}g_{i_1,\ldots,i_n}(x_1^p,x_2^p,\ldots,x_n^p)+g_0(x_1^p,x_2^p,\ldots,x_n^p),$$
where $g_{i_1,\ldots,i_n},~g_0\in K[x_1^p,\ldots,x_n^p]$.

$(1)$ If $f_n\neq 0$, then $D(x_n)=f_n\neq 0$. Thus, we have $1\in \operatorname{Im}D$. Therefore, we have that
$$D(g)=\sum_{q=1}^nf_q(x_{q+1},\ldots,x_n)(\sum_{\mbox {\tiny $\begin{array}{c}
0\leq i_1\leq p-1\\
\vdots\\
0\leq i_n\leq p-1\\
i_1+\cdots+i_n\geq 1\end{array}$}}i_qx_1^{i_1}\cdots x_q^{i_q-1}\cdots x_n^{i_n}g_{i_1,\ldots,i_n}(x_1^p,\ldots,x_n^p)).$$
Since $0\leq i_q\leq p-1$ for all $1\leq q\leq n$ and $i_1+\cdots+i_n\geq 1$, we have that the monomial $x_1^{p-1}x_2^{p-1}\cdots x_n^{p-1}$ doesn't appear in $D(g)$ for any $g\in K[x]$. Hence $x_1^{p-1}x_2^{p-1}\cdots x_n^{p-1}\notin \operatorname{Im}D$, whence $\operatorname{Im}D$ is not a Mathieu-Zhao space of $K[x]$.

$(2)$ If $f_n=0$, then we can assume that $f_n=f_{n-1}=\cdots=f_{n-j}=0$ for some $j\in \{0,1,\ldots,n-2\}$ and $f_{n-j-1}\neq 0$. Then we have
$$D(g)=\sum_{q=1}^{n-j-1}f_q(x_{q+1},\ldots,x_n)(\sum_{\mbox {\tiny $\begin{array}{c}
0\leq i_1\leq p-1\\
\vdots\\
0\leq i_n\leq p-1\\
i_1+\cdots+i_n\geq 1\end{array}$}}i_qx_1^{i_1}\cdots x_q^{i_q-1}\cdots x_n^{i_n}g_{i_1,\ldots,i_n}(x_1^p,\ldots,x_n^p)).$$
Since $D(x_{n-j-1}f_{n-j-1}^i)=f_{n-j-1}^{i+1}$
for all $i\in \mathbb{N}$, we have $f_{n-j-1}^m\in \operatorname{Im}D$ for all $m\geq 1$. Since $0\leq i_q\leq p-1$ for all $1\leq q\leq n-j-1$ and $D(g_{0,\ldots,0,i_{n-j},\ldots,i_n}(x_1^p,\ldots,x_n^p))=0$, the polynomial $x_1^{p-1}x_2^{p-1}\cdots x_{n-j-1}^{p-1}f_{n-j-1}^m$ doesn't appear in $D(g)$ for any $g\in K[x]$, for all $m\geq 1$. Hence $x_1^{p-1}x_2^{p-1}\cdots x_{n-j-1}^{p-1}f_{n-j-1}^m\notin \operatorname{Im}D$ for all $m\geq 1$, whence $\operatorname{Im}D$ is not a Mathieu-Zhao space of $K[x]$. Therefore, we have $D=0$.
\end{proof}

\section{Problem \ref{prob1.5} for the case $n=1$}

\begin{prop}\label{prop3.1}
Let $D=x_1^rf_1(x_1^p)\partial_1$ be a nonzero derivation of $R[x_1]$ and $f_1(x_1^p)$ a polynomial of $x_1^p$. Then $D$ is locally nilpotent if and only if $r\neq lp+1$ for some $l\in \mathbb{N}$.
\end{prop}
\begin{proof}
$``\Rightarrow"$ It suffices to prove that if $r=lp+1$ for some $l\in \mathbb{N}$, then $D$ is not locally nilpotent. If $r=lp+1$, then $D=x_1^{lp+1}f_1(x_1^p)\partial_1=x_1\tilde{f}_1(x_1^p)\partial_1$, where $\tilde{f}_1(x_1^p)=x_1^{lp}f_1(x_1^p)$. Thus, we have
$$D^M(x_1)=\tilde{f}_1(x_1^p)^Mx_1$$
for any $M\in \mathbb{N}^*$. Hence $D$ is not locally nilpotent.

$``\Leftarrow"$ It suffices to prove that if $r=0,2,\ldots,p-1$, then $D$ is locally nilpotent. If $r=0$, then $D=f_1(x_1^p)\partial_1$. Thus, we have $D^{i+1}(x_1^i)=0$ for any $i\in \mathbb{N}$. Thus, $D$ is locally nilpotent. If $r\in \{2,3,\ldots,p-1\}$, then $D^j(x_1)=[(j-1)r-j+2][(j-2)r-j+3]\cdots (2r-1)rx_1^{jr-(j-1)}f_1(x_1^p)^j$. Thus, we have
$$D^j(x_1)=[(j-1)(r-1)+1][(j-2)(r-1)+1]\cdots[(r-1)+1]x_1^{jr-(j-1)}f_1(x_1^p)^j.$$
Since $(r-1,p)=1$, it follows from Chinese remainder theorem that there exists $N_r\in \mathbb{N}$ such that $N_r=1\operatorname{mod}(r-1)$ and $N_r=0\operatorname{mod}p$ for all $2\leq r\leq p-1$. Let $J=\operatorname{max}\{N_2,N_3,\ldots,N_{p-1}\}+1$. Then we have $D^J(x_1)=0$ for all $2\leq r\leq p-1$. Thus, $D$ is locally nilpotent.
\end{proof}

\begin{cor}\label{cor3.2}
Let $D=x_1^rf_1(x_1^p)\partial_1$ be a nonzero derivation of $R[x_1]$, where $f_1(x_1^p)$ is a polynomial of $x_1^p$. Then $\operatorname{Im}D$ is a Mathieu-Zhao space of $R[x_1]$ if and only if $D$ is not locally nilpotent.
\end{cor}
\begin{proof}
The conclusion follows from Proposition \ref{prop2.3} and Proposition \ref{prop3.1}.
\end{proof}

\begin{rem}\label{rem3.3}
If $r=1$ and $f_1(x_1^p)\in R^*$ in Proposition \ref{prop3.1}, then $D$ is locally finite. If $r= 1$ and $f_1(x_1^p)\in R[x_1^p]/R$, then $D$ is not locally finite.
\end{rem}

\begin{thm}\label{thm3.4}
Let $R$ be an integral domain of characteristic 2 and $D=f(x_1)\partial_1$ a nonzero derivation of $R[x_1]$. Then $D$ is locally nilpotent if and only if $f(x_1)=f_0(x_1^2)$ for some $f_0(x_1^2)\in R[x_1^2]$.
\end{thm}
\begin{proof}
$``\Leftarrow"$ If $f(x_1)=f_0(x_1^2)$ for some $f_0(x_1^2)\in R[x_1^2]$, then $D=f_0(x_1^2)\partial_1$. Thus, we have $D^2(x_1)=0$. Hence $D$ is locally nilpotent.

$``\Rightarrow"$ Since $f(x_1)\in R[x_1]$, we have $f(x_1)=f_0(x_1^2)$ or $x_1f_1(x_1^2)$ or $a_0f_0(x_1^2)+a_1x_1f_1(x_1^2)$ for some $a_0,a_1\in R^*$, $f_0,f_1\in R[x_1^2]$. It suffices to prove that if $f(x_1)=x_1f_1(x_1^2)$ or $a_0f_0(x_1^2)+a_1x_1f_1(x_1^2)$, then $D$ is not locally nilpotent.

$(1)$ If $f(x_1)=x_1f_1(x_1^2)$, then $D=x_1f_1(x_1^2)\partial_1$ and $D^M(x_1)=x_1f_1(x_1^2)^M$ for any $M\in \mathbb{N}^*$. Thus, $D$ is not locally nilpotent.

$(2)$ If $f(x_1)=a_0f_0(x_1^2)+a_1x_1f_1(x_1^2)$, then $D=(a_0f_0(x_1^2)+a_1x_1f_1(x_1^2))\partial_1$ and $D^M(x_1)=[a_0f_0(x_1^2)+a_1x_1f_1(x_1^2)](a_1f_1(x_1^2))^{M-1}$
for any $M\in \mathbb{N}$, $M\geq 2$. Thus, $D$ is not locally nilpotent. Then the conclusion follows.
 \end{proof}

\begin{cor}\label{cor3.5}
Let $R$ be an integral domain of characteristic 2 and $D=f(x_1)\partial_1$ a nonzero derivation of $R[x_1]$. Then $\operatorname{Im}D$ is a Mathieu-Zhao space of $R[x_1]$ if and only if $D$ is not locally nilpotent.
\end{cor}
\begin{proof}
The conclusion follows from Theorem \ref{thm2.4} and Theorem \ref{thm3.4}.
\end{proof}

\begin{thm}\label{thm3.6}
Let $R$ be an integral domain of characteristic $\leq 3$ and $D=f(x_1)\partial_1$ a nonzero derivation of $R[x_1]$. Then $D$ is locally nilpotent if and only if $f(x_1)=x_1^r\tilde{g}(x_1^p)$ for some $r\in \{0,2,\ldots,p-1\}$, $\tilde{g}(x_1^p)\in R[x_1^p]$ or $f(x_1)=x_1^2f_2(x_1^3)+x_1f_1(x_1^3)+f_0(x_1^3)$ and $f_2f_0=f_1^2$ for $f_2,f_1,f_0\in R[x_1^3]$ and $f_2f_1f_0\neq 0$.
\end{thm}
\begin{proof}
If $p=2$, then the conclusion follows from the proof of Theorem \ref{thm3.4}.

If $p=3$, then $f(x_1)$ can be written as $f_0(x_1^3)$, $x_1f_1(x_1^3)$, $x_1^2f_1(x_1^3)$, $c_2x_1^2f_2(x_1^3)+c_1x_1f_1(x_1^3)$, $c_2x_1^2f_2(x_1^3)+c_0f_0(x_1^3)$, $c_1x_1f_1(x_1^3)+c_0f_0(x_1^3)$ or $c_2x_1^2f_2(x_1^3)+c_1x_1f_1(x_1^3)+c_0f_0(x_1^3)$ with $c_0, c_1, c_2\in R^*$, $f_0, f_1, f_2\in R[x_1^3]$ and $f_0f_1f_2\neq 0$. Without loss of generality, we can assume that $c_0=c_1=c_2=1$.

If $f(x_1)=f_0(x_1^3)$, $x_1f_1(x_1^3)$, $x_1^2f_1(x_1^3)$, then the conclusion follows from Proposition \ref{prop2.3}.

If $f(x_1)=x_1^2f_2(x_1^3)+x_1f_1(x_1^3)$, then we claim that $D^{2k}(x_1)=(2x_1^3f_2^2+x_1f_1^2)f_1^{2(k-1)}$ and $D^{2k+1}(x_1)=(x_1^2f_2+x_1f_1)f_1^{2k}$ for all $k\in \mathbb{N}^*$. It's easy to check that if $k=1$, then the conclusion follows. Suppose that the conclusion is true for $k=l$. If $k=l+1$, then $D^{2l+2}(x_1)=D(D^{2l+1}(x_1))=(2x_1^3f_2^2+x_1f_1^2)f_1^{2l}$
and $D^{2l+3}(x_1)=D(D^{2l+2}(x_1))=(x_1^2f_2+x_1f_1)f_1^{2(l+1)}.$
Thus, $D$ is not locally nilpotent.

If $f(x_1)=x_1^2f_2(x_1^3)+f_0(x_1^3)$, then we claim that $D^{2k}(x_1)=2^k(x_1^3f_2+x_1f_0)f_0^{k-1}f_2^k$ and $D^{2k+1}(x_1)=2^k(x_1^2f_2+f_0)(f_0f_2)^k$ for all $k\in \mathbb{N}^*$. If $k=1$, then it's easy to check that the conclusion follows. Suppose that the conclusion is true for $k=l$. If $k=l+1$, then $D^{2l+2}(x_1)=D(D^{2l+1}(x_1))=2^{l+1}(x_1^3f_2+x_1f_0)f_0^lf_2^{l+1}$ and $D^{2l+3}(x_1)=D(D^{2l+2}(x_1))=2^{l+1}(x_1^2f_2+f_0)(f_0f_2)^{l+1}$. Thus, $D$ is not locally nilpotent.

If $f(x_1)=x_1f_1(x_1^3)+f_0(x_1^3)$, then we claim that $D^k(x_1)=(x_1f_1+f_0)f_1^{k-1}$ for all $k\in \mathbb{N}^*$. If $k=1$, then it's easy to check that the conclusion is true. Suppose that the conclusion is true for $k=l$. If $k=l+1$, then $D^{l+1}(x_1)=D(D^l(x_1))=(x_1f_1+f_0)f_1^l$. Thus, $D$ is not locally nilpotent.

If $f(x_1)=x_1^2f_2(x_1^3)+x_1f_1(x_1^3)+f_0(x_1^3)$, then we claim that $D^{2k}(x_1)=(x_1^2f_2+x_1f_1+f_0)(2x_1f_2+f_1)(2f_0f_2+f_1^2)^{k-1}$ and $D^{2k+1}(x_1)=(x_1^2f_2+x_1f_1+f_0)(2f_0f_2+f_1^2)^k$ for all $k\in \mathbb{N}^*$. If $k=1$, then it's easy to check that the conclusion follows. Suppose that the conclusion is true for $k=l$. If $k=l+1$, then $D^{2l+2}(x_1)=D(D^{2l+1}(x_1))=(x_1^2f_2+x_1f_1+f_0)(2x_1f_2+f_1)(2f_0f_2+f_1^2)^l$ and $D^{2l+3}(x_1)=D(D^{2l+2}(x_1))=(x_1^2f_2+x_1f_1+f_0)(2f_0f_2+f_1^2)^{l+1}$.
If $f_1^2=f_0f_2$, then $D^3(x_1)=0$. Thus, $D$ is locally nilpotent. If $f_1^2\neq f_0f_2$, then $D$ is not locally nilpotent. Hence the conclusion follows.
\end{proof}

\begin{lem}\label{lem3.7}
Let $D=f(x_1)\partial_1$ be a nonzero derivation of $R[x_1]$. If $f(x_1)=c_{i_1}x_1^{i_1}f_{i_1}(x_1^p)+c_{i_2}x_1^{i_2}f_{i_2}(x_1^p)$ for $c_{i_1}, c_{i_2}\in R^*$, $f_{i_1}, f_{i_2}\in R[x_1^p]$, $f_{i_1}f_{i_2}\neq 0$, $i_1\neq i_2$, $i_1, i_2\in \{0,1,\ldots,p-1\}$, then $$D^ {k+1}(x_1)=(ki_1-(k-1))a_{k-1,1}x_1^{(k+1)i_1-k}f_{i_1}^{k+1}+
\sum_{j=2}^{k+1}((k-j+1)i_1+(j-1)i_2-(k-1))\cdot$$
$$a_{k-1,j}x_1^{(k-j+2)i_1+(j-1)i_2-k}f_{i_1}^{k-j+2}f_{i_2}^{j-1}
+(ki_2-(k-1))a_{k-1,k+1}x_1^{(k+1)i_2-k}f_{i_2}^{k+1}$$
for all $k\in \mathbb{N}^*$, where $a_{k-1,1}, a_{k-1,2},\ldots,a_{k-1,k+1}$ are the coefficients of $x_1^{ki_1-(k-1)}f_{i_1}^k$, $x_1^{(k-1)i_1+i_2-(k-1)}f_{i_1}^{k-1}f_{i_2}$,\ldots,$x_1^{ki_2-(k-1)}f_{i_2}^{k+1}$
in $D^k(x_1)$ respectively and $a_{0,1}=c_{i_1}$, $a_{0,2}=c_{i_2}$.
\end{lem}
\begin{proof}
Since $D(x_1)=c_{i_1}x_1^{i_1}f_{i_1}(x_1^p)+c_{i_2}x_1^{i_2}f_{i_2}(x_1^p)$, we have $a_{0,1}=c_{i_1}$, $a_{0,2}=c_{i_2}$. Without loss of generality, we can assume that $c_{i_1}=c_{i_2}=1$. Since $D^2(x_1)=i_1x_1^{2i_1-1}f_{i_1}^2+(i_1+i_2)x_1^{i_1+i_2-1}f_{i_1}f_{i_2}+i_2x_1^{2i_2-1}f_{i_2}^2$, we have $a_{1,1}=i_1a_{0,1}$, $a_{1,2}=i_1a_{0,2}+i_2a_{0,1}$, $a_{1,3}=i_2a_{0,2}$. Thus, the conclusion is true for $k=1$. Suppose that the conclusion is true for $k=l$. If $k=l+1$, then
$$D^{l+2}(x_1)=D(D^{l+1}(x_1))=(x_1^{i_1}f_{i_1}+x_1^{i_2}f_{i_2})[((l+1)i_1-l)a_{l,1}x_1^{(l+1)i_1-l-1}f_{i_1}^{l+1}
+$$
$$(li_1+i_2-l)a_{l,2}x_1^{li_1+i_2-l-1}f_{i_1}^lf_{i_2}+\cdots+(i_1+li_2-l)a_{l,l+1}x_1^{i_1+li_2-l-1}f_{i_1}f_{i_2}^l+((l+1)i_2-l)\times$$
$$a_{l,l+2}x_1^{(l+1)i_2-l-1}f_{i_2}^{l+1}]
=[(l+1)i_1-l]a_{l,1}x_1^{(l+2)i_1-l-1}f_{i_1}^{l+2}
+[(li_1+i_2-l)a_{l,2}+((l+1)i_1-$$
$$l)a_{l,1}]x_1^{(l+1)i_1+i_2-l-1}f_{i_1}^{l+1}f_{i_2}+\cdots+[((l+1)i_2-l)a_{l,l+2}+(i_1+li_2-l)a_{l,l+1}]x_1^{i_1+(l+1)i_2-l-1}\times$$
$$f_{i_1}f_{i_2}^{l+1}+[(l+1)i_2-l]a_{l,l+2}x_1^{(l+2)i_2-l-1}f_{i_2}^{l+2}.$$
Thus, we have $a_{l+1,1}=((l+1)i_1-l)a_{l,1}$, $a_{l+1,j}=((l-j+2)i_1+(j-1)i_2-l)a_{l,j}+((l-j+3)i_1+(j-2)i_2-l)a_{l,j-1}$ and $a_{l+1,l+3}=((l+1)i_2-l)a_{l,l+2}$ for $2\leq j\leq l+2$, which completes the proof.
\end{proof}

\begin{prop}\label{prop3.8}
Let $R$ and $D$ be as in Lemma \ref{lem3.7}. If $i_1=1$ or $i_2=1$, then $D$ is not locally nilpotent.
\end{prop}
\begin{proof}
Let $a_{k,1}, a_{k,2},\ldots,a_{k,k+2}$ be the coefficients of $x_1^{(k+1)i_1-k}f_{i_1}^{k+1}$, $x_1^{ki_1+i_2-k}f_{i_1}^kf_{i_2}$,\allowbreak\ldots,$x_1^{(k+1)i_2-k}f_{i_2}^{k+1}$
in $D^{k+1}(x_1)$ respectively.

If $i_1=1$, then $a_{k,1}=(ki_1-(k-1))a_{k-1,1}=a_{k-1,1}$ for all $k\in \mathbb{N}^*$. Since $a_{0,1}=c_{i_1}$, we have $a_{k,1}=c_{i_1}\in R^*$ for all $k\in \mathbb{N}^*$. Thus, $D$ is not locally nilpotent.

If $i_2=1$, then $a_{k,k+2}=(ki_2-(k-1))a_{k-1,k+1}=a_{k-1,k+1}$ for all $k\in \mathbb{N}^*$. Since $a_{0,2}=c_{i_2}$, we have $a_{k,k+2}=c_{i_2}\in R^*$ for all $k\in \mathbb{N}^*$. Thus, $D$ is not locally nilpotent.
\end{proof}

\section{Problem \ref{prob1.4} for derivations of $K[x_1]$}

\begin{prop}\label{prop4.1}
Let $D=x_1^if_1(x_1^p)\partial_1$ be a derivation of $K[x_1]$ with $f_1(x_1^p)\in K[x_1^p]$ for $0\leq i\leq p-1$. Suppose that $I$ is a nonzero ideal generated by $x_1^jg_1(x_1^p)$ for $0\leq j\leq p-1$, $g_1(x_1^p)\in K[x_1^p]$. Then we have the following statements:

$(1)$ If $i=1$, then $\mathfrak{r}(D(I))=0$. In particular, $D(I)$ is a Mathieu-Zhao space of $K[x_1]$.

$(2)$ If $i\neq 1$, then $D(I)$ is not a Mathieu-Zhao space of $K[x_1]$.
\end{prop}
\begin{proof}
Since $I=(x_1^jg_1(x_1^p))$, we have that
$$I=K\{x_1^{j+kp}g_1(x_1^p),~ x_1^{j+kp+1}g_1(x_1^p),\ldots,x_1^{j+kp+p-1}g_1(x_1^p)~\operatorname{for}~\operatorname{all}~k\in \mathbb{N}\}.$$
Then we have
$$D(I)=\{jx_1^{i+j+kp-1}g_1f_1,~ (j+1)x_1^{i+j+kp}g_1f_1,\ldots,(j+p-1)x_1^{i+j+kp+p-2}g_1f_1~\operatorname{for}~\operatorname{all}$$ $$
k\in\mathbb{N}\}.$$
If $j+t-1=p-1$, then $j+t=0$ for any $t\in\{0,1,\ldots,p-2\}$. Thus, we have $x_1^{i+kp+p-1}g_1f_1\notin D(I)$ for all $k\in\mathbb{N}$ and $x_1^{i+kp+l}g_1f_1\in D(I)$ for all $k\in \mathbb{N}$ for any $l\in\{0,1,\ldots,p-2\}$.

$(1)$ If $i=1$, then $x_1^{kp}g_1f_1\notin D(I)$ for all $k\in\mathbb{N}$. We claim that $\mathfrak{r}(D(I))=0$. Suppose that there exists nonzero polynomial $\alpha(x_1)\in \mathfrak{r}(D(I))$. Then we have that $\alpha^m\in D(I)$ for all $m\geq N$. Let $Q(x_1):=\alpha^N$. Then $Q(x_1)^m\in D(I)$ for all $m\in \mathbb{N}^*$ and $Q(x_1)$ can be written as
$$Q(x_1)=(b_{t_1}x_1^{i+j+t_1-1}+b_{t_2}x_1^{i+j+t_2-1}+\cdots+b_{t_q}x_1^{i+j+t_q-1})x_1^{kp}g_1f_1$$
for some $k\in \mathbb{N}$, $0\leq t_1\leq t_2\leq\cdots\leq t_q\leq p-1$, $b_{t_1}, b_{t_2},\ldots,b_{t_q}\in K$. Hence we have $Q(x_1)^{p^m}\notin D(I)$ for all $m\in \mathbb{N}^*$. This is a contradiction. Thus, we have $\alpha(x_1)=0$. That is, $\mathfrak{r}(D(I))=0$.

$(2)$ If $i\neq 1$, then $x_1^{kp}g_1f_1\in D(I)$ for all $k\geq 2$. Since $(x_1^{2p}g_1f_1)^m=(x_1^{2p}g_1f_1)^{m-1}\allowbreak (x_1^{2p}g_1f_1)$ and $(x_1^{2p}g_1f_1)^{m-1}$ is a polynomial of $x_1^p$, we have $(x_1^{2p}g_1f_1)^m\in D(I)$ for all $m\in \mathbb{N}^*$. It's easy to see that $x_1^{i+p-1}(x_1^{2p}g_1f_1)^m\notin D(I)$ for all $m\in \mathbb{N}^*$. Hence $D(I)$ is not a Mathieu-Zhao space of $K[x_1]$.
\end{proof}

\begin{prop}\label{prop4.2}
Let $D=x_1^{i_1}f_1(x_1^p)\partial_1$ be a derivation of $K[x_1]$ with $f_1(x_1^p)\in K[x_1^p]$ for $i_1\in \{0,1,\ldots,p-1\}$, $f_1(x_1^p)\in K[x_1^p]$ and $f_1(x_1^p)\neq 0$. Suppose that $I$ is any nonzero ideal of $K[x_1]$. Then we have the following statements:

$(1)$ If $I$ is generated by $x_1^{j_1}g_1(x_1^p)$ for $0\leq j_1\leq p-1$, $g_1(x_1^p)\in K[x_1^p]$, then $D(I)$ is a Mathieu-Zhao space of $K[x_1]$ for $i_1=1$
and $D(I)$ is not a Mathieu-Zhao space of $K[x_1]$ for $i_1\neq 1$.

$(2)$ If $I$ is generated by $x_1^{j_1}g_1(x_1^p)+x_1^{j_2}g_2(x_1^p)$ for $0\leq j_1<j_2\leq p-1$, $g_1, g_2\in K[x_1^p]$, then $D(I)$ is a Mathieu-Zhao space of $K[x_1]$ for $j_2-j_1\neq i_1-1$ and $p-i_1+1$ and $D(I)$ is not a Mathieu-Zhao space of $K[x_1]$ for $j_2-j_1=i_1-1$ or $p-i_1+1$.

$(3)$ If $I$ is generated by $x_1^{j_1}g_1(x_1^p)+x_1^{j_2}g_2(x_1^p)+\cdots+x_1^{j_r}g_r(x_1^p)$ for $0\leq j_1<j_2<\cdots<j_r\leq p-1$, $g_1, g_2,\ldots,g_r\in K[x_1^p]$ and $r\geq 3$, then $D(I)$ is a Mathieu-Zhao space of $K[x_1]$.
\end{prop}
\begin{proof}
$(1)$ The conclusion follows from Proposition \ref{prop4.1}.

$(2)$ It's easy to compute that
$$D(I)=K\{x_1^{kp}[(t+j_1)x_1^{t+j_1+i_1-1}g_1f_1+(t+j_2)x_1^{t+j_2+i_1-1}g_2f_1]|
0\leq t\leq p-1,~k\in\mathbb{N}\}.$$
Suppose that $x_1^{kp}g_if_1\in D(I)$ for some $k\in \mathbb{N}$, $i=1$ or 2. Then we have the following equations:
\begin{equation}
  \left\{ \begin{aligned}
  \nonumber
  t+j_1=0~\operatorname{or}~p~~~~~~~~~~~~~~~~~~~~~~~~~~~~~~~~~~~~~~~~~~~~\\
  t+j_2+i_1-1=0~\operatorname{or}~p~\operatorname{or}~2p~~~~~~~~~~~~~~~~~~~~~~~~~\\
                          \end{aligned} \right.
  \end{equation}
or
\begin{equation}
  \left\{ \begin{aligned}
  \nonumber
  t+j_1+i_1-1=0~\operatorname{or}~p~\operatorname{or}~2p~~~~~~~~~~~~~~~~~~~~~~~~~\\
  t+j_2=0~\operatorname{or}~p~~~~~~~~~~~~~~~~~~~~~~~~~~~~~~~~~~~~~~~~~~~~\\
                          \end{aligned} \right.
  \end{equation}
Since $0\leq t, j_1, j_2, i_1\leq p-1$, the above equations are equivalent to that $j_2-j_1=i_1-1$ or $p-i_1+1$. Since $j_2>j_1$, we have $i_1\neq 1$ in this situation.

$(2.1)$ If $j_2-j_1=i_1-1$, then
$$D(I)=K\{x_1^{kp}[(t+j_1)x_1^{t+j_1+i_1-1}g_1f_1+(t+j_1+i_1-1)x_1^{t+j_1+2i_1-2}g_2f_1]|~0\leq t\leq p-1,$$
$$k\in \mathbb{N}\}.$$
Since $0\leq t\leq p-1$, there exists $t\in\{0,1,\ldots,p-1\}$ such that $p|(t+j_1+i_1-1)$ and $t+j_1\neq 0$. Thus, we have $x_1^{(k+3)p}g_1f_1\in D(I)$ for all $k\in \mathbb{N}$. Since $(x_1^{3p}g_1f_1)^m=(x_1^{3p}g_1f_1)^{m-1}\cdot(x_1^{3p}g_1f_1)$ and $(x_1^{3p}g_1f_1)^{m-1}$ is a polynomial of $x_1^p$, we have $(x_1^{3p}g_1f_1)^m\in D(I)$ for all $m\in \mathbb{N}^*$. It's easy to check that $x_1^{p-1}(x_1^{3p}g_1f_1)^m\notin D(I)$ for all $m\in \mathbb{N}^*$. Thus, $D(I)$ is not a Mathieu-Zhao space of $K[x_1]$.

$(2.2)$ If $j_2-j_1=p-i_1+1$, then
$$D(I)=K\{x_1^{kp}[(t+j_1)x_1^{t+j_1+i_1-1}g_1f_1+(t+j_1-i_1+1)x_1^{t+p+j_1}g_2f_1]|0\leq t\leq p-1,~k\in \mathbb{N}\}.$$
Since $0\leq t\leq p-1$, there exists $t\in\{0,1,\ldots,p-1\}$ such that $p|(t+j_1)$ and $t+j_1-i_1+1 \neq 0$. Thus, we have $x_1^{(k+3)p}g_2f_1\in D(I)$ for all $k\in \mathbb{N}$. Since $(x_1^{3p}g_2f_1)^m=(x_1^{3p}g_2f_1)^{m-1}\cdot(x_1^{3p}g_2f_1)$ and $(x_1^{3p}g_2f_1)^{m-1}$ is a polynomial of $x_1^p$, we have $(x_1^{3p}g_2f_1)^m\in D(I)$ for all $m\in \mathbb{N}^*$. It's easy to check that $x_1^{p-1}(x_1^{3p}g_2f_1)^m\notin D(I)$ for all $m\in \mathbb{N}^*$. Thus, $D(I)$ is not a Mathieu-Zhao space of $K[x_1]$.

$(2.3)$ If $j_2-j_1\neq i_1-1$ and $p-i_1+1$, then $x_1^{kp}g_1f_1, x_1^{kp}g_2f_1\notin D(I)$ for any $k\in \mathbb{N}$. It's easy to check that $\tilde{h} \neq \tilde{c}\operatorname{mod} x_1^p$ for any nonzero $\tilde{h}\in D(I)$, $\tilde{c}\in K$. We claim that $\mathfrak{r}(D(I))=0$. Assume otherwise, if there exists a nonzero polynomial $q(x_1)$ such that $q(x_1)\in \mathfrak{r}(D(I))$, then $q(x_1)$ can be written as
$$q(x_1)=x_1^{k_1}g_{k_1}(x_1^p)+\cdots+x_1^{k_w}g_{k_w}(x_1^p)$$
for $w\geq 1$, $g_{k_1},\ldots,g_{k_w}\in K[x_1^p]$, $g_{k_1}\cdots g_{k_w}\neq 0$, $k_1,\ldots,k_w$ are distinct integers between 0 and $p-1$. Thus, we have that
$$q(x_1)^{p^m}=x_1^{k_1p^m}g_{k_1}(x_1^p)^{p^m}+\cdots+x_1^{k_wp^m}g_{k_w}(x_1^p)^{p^m}=\bar{c}\operatorname{mod} x_1^p$$
for some $\bar{c}\in K$ and for all $m\in \mathbb{N}^*$. Hence we have $q(x_1)^{p^m}\notin D(I)$ for all $m\in \mathbb{N}^*$, whence $q(x_1)\notin \mathfrak{r}(D(I))$, which is a contradiction. Therefore, we have that $\mathfrak{r}(D(I))=0$. Hence $D(I)$ is a Mathieu-Zhao space of $K[x_1]$.

$(3)$ It's easy to check that
$$D(I)=K\{x_1^{kp}[(t+j_1)x_1^{t+i_1+j_1-1}g_1f_1+\cdots+(t+j_r)x_1^{t+i_1+j_r-1}g_rf_1]|0\leq t\leq p-1,~k\in \mathbb{N}\}.$$
Since $r\geq 3$, $\tilde{h}$ contains at least two different monomials by module $x_1^p$ for any nonzero $\tilde{h}\in D(I)$. We have that $\mathfrak{r}(D(I))=0$ by following the arguments of $(2.3)$ of Proposition \ref{prop4.2}. Thus, $D(I)$ is a Mathieu-Zhao space of $K[x_1]$.
\end{proof}

\begin{cor}\label{cor4.3}
Let $D=c\partial_1$ be a derivation of $K[x_1]$ with $c\in K^*$. Suppose that $I$ is any nonzero ideal of $K[x_1]$. Then we have the following statements:

$(1)$ If $I$ is generated by $x_1^lg_l(x_1^p)$ or $a_{i_1}x_1^{i_1}g_{i_1}(x_1^p)+a_{i_1+1}x_1^{i_1+1}g_{i_1+1}(x_1^p)$ or $a_0g_0(x_1^p)\allowbreak+a_{p-1}x_1^{p-1}g_{p-1}(x_1^p)$ for some $l\in\{0,1,\ldots,p-1\}$,
$a_0, a_{i_1}, a_{i_1+1}, a_{p-1}\in K^*$, $g_l(x_1^p)\in K[x_1^p]$ and $g_l(x_1^p)\neq 0$ for all $l\in\{0,1,\ldots,p-1\}$, $i_1\in\{0,1,\ldots,p-2\}$, then
$D(I)$ is not a Mathieu-Zhao space of $K[x_1]$. In particular, if $p=2$, then $D(I)$
is not a Mathieu-Zhao space of $K[x_1]$ for any ideal $I$ of $K[x_1]$.

$(2)$ If $I$ is not the cases of $(1)$, then $\mathfrak{r}(D(I))=0$. In particular, $D(I)$ is a Mathieu-Zhao space of $K[x_1]$.
\end{cor}
\begin{proof}
The conclusion follows from the proof of Proposition \ref{prop4.2}.
\end{proof}

\begin{thm}\label{thm4.4}
Let $D=f(x_1)\partial_1$ be a nonzero derivation of $K[x_1]$ with $f(x_1)\in K[x_1]$. Suppose that $I$ is a nonzero ideal of $K[x_1]$. Then we have the following statements:

$(1)$ If $f(x_1)=x_1^{i_1}f_1(x_1^p)$ for some $i_1\in \{0,1,\ldots,p-1\}$, $f_1\in K[x_1^p]$, then the conclusion is the same as in Proposition \ref{prop4.2}.

$(2)$ If $f(x_1)=x_1^{i_1}f_1(x_1^p)+\cdots+x_1^{i_s}f_s(x_1^p)$ for $s\geq 2$, $0\leq i_1<\cdots<i_s\leq p-1$, $f_1,\ldots,f_s\in K[x_1^p]$, then $\mathfrak{r}(D(I))=0$. In particular, $D(I)$ is a Mathieu-Zhao space of $K[x_1]$.
\end{thm}
\begin{proof}
$(1)$ The conclusion follows from Proposition \ref{prop4.2}.

$(2)$ Since $K[x_1]$ is a principal ideal domain, we can assume that $I$ is generated by $x_1^{j_1}g_1(x_1^p)+x_1^{j_2}g_2(x_1^p)+\cdots+x_1^{j_r}g_r(x_1^p)$ for $r\geq 1$, $0\leq j_1<j_2<\cdots<j_r\leq p-1$, $g_1,\ldots,g_r\in K[x_1^p]$. Then we have that
$D(I)=K\{x_1^{kp}[(t+j_1)(x_1^{t+j_1+i_1-1}f_1g_1+\cdots+x_1^{t+j_1+i_s-1}f_sg_1)+\cdots+
(t+j_r)(x_1^{t+j_r+i_1-1}f_1g_r+\cdots+x_1^{t+j_r+i_s-1}f_sg_r)]|0\leq t\leq p-1,~k\in \mathbb{N}\}$. Thus, we have that $\tilde{h}$ contains at least two different monomials by module $x_1^p$ for any nonzero $\tilde{h}\in D(I)$. We deduce that $\mathfrak{r}(D(I))=0$ by following the arguments of $(3)$ of Proposition \ref{prop4.2}. Therefore, $D(I)$ is a Mathieu-Zhao space of $K[x_1]$.
\end{proof}

\begin{cor}\label{cor4.5}
Let $K$ be a field of characteristic 2 and $D=f(x_1)\partial_1$ a nonzero derivation of $K[x_1]$. Then $D(I)$ is a Mathieu-Zhao space of $K[x_1]$ if and only if $f(x_1)\neq f_0(x_1^2)$ for any $f_0(x_1^2)\in K[x_1^2]$, where $I$ is any nonzero ideal of $K[x_1]$.
\end{cor}
\begin{proof}
The conclusion follows from Theorem \ref{thm4.4}.
\end{proof}

\begin{cor}\label{cor4.6}
Let $K$ be a field of characteristic 2 and $D=f(x_1)\partial_1$ a nonzero derivation of $K[x_1]$. Then $D(I)$ is a Mathieu-Zhao space of $K[x_1]$ if and only if $D$ is not locally nilpotent, where $I$ is any nonzero ideal of $K[x_1]$.
\end{cor}
\begin{proof}
The conclusion follows from Theorem \ref{thm3.4} and Corollary \ref{cor4.5}.
\end{proof}

Next we give partial answer to Problem \ref{prob1.4} for $\mathcal{E}$-derivations $\delta$ of $K[x_1]$, where $\delta=I-\phi$, $I$ is the identity map. Hence we use the notation $I_0$ for any ideal of $K[x_1]$ in the following discussion.

\begin{thm}\label{thm4.7}
Let $\delta=I-\phi$ be a $\mathcal{E}$-derivation of $K[x_1]$ and $I_0$ a nonzero ideal of $K[x_1]$. Then we have the following statements:

$(1)$ If $I_0$ is generated by $x_1+a$ for some $a\in K$, then $\delta(I_0)$ is not a Mathieu-Zhao space of $K[x_1]$ in the case that $\phi(x_1)=x_1+c$ for some $c\in K^*$. Otherwise, $\delta(I_0)$ is a Mathieu-Zhao space of $K[x_1]$.

$(2)$ If $I_0$ is generated by $x_1^i$ for some integer $i\geq 2$, then $\delta(I_0)$ is not a Mathieu-Zhao space of $K[x_1]$ in the case that $\phi(x_1)=c$ or $\phi(x_1)=x_1+c$ for some $c\in K^*$. Otherwise, $\delta(I_0)$ is a Mathieu-Zhao space of $K[x_1]$.
\end{thm}
\begin{proof}
It follows from the discussion after Lemma 3.1 in \cite{4} that $\delta=I-\phi$ can be divided into four(exhausting) cases: $\phi(x_1)=c$ for some $c\in K$ or $\phi(x_1)=x_1+c$ for some $c\in K$ or $\phi(x_1)=\tilde{q}x_1$ for some $\tilde{q}\in K$ or $\operatorname{deg}\phi(x_1)\geq 2$.

If $\operatorname{deg}\phi(x_1)\geq 2$, then it follows from Proposition 3.7 in \cite{4} that $\delta(I_0)$ is a Mathieu-Zhao space of $K[x_1]$ for any ideal $I_0$ of $K[x_1]$.\\

$(1)$ $(1.1)$ If $\phi(x_1)=c$, then $\delta(I_0)$ is an ideal generated by $x_1-c$. Clearly, $\delta(I_0)$ is a Mathieu-Zhao space of $K[x_1]$.

$(1.2)$ If $\phi(x_1)=x_1+c$ and $c\neq 0$, then $\delta(x_1+a)=-c$. Thus, we have $1\in \delta(I_0)$. We claim that $x_1^{p-1}\notin \delta(I_0)$. Suppose that $x_1^{p-1}\in \delta(I_0)$. Then there exists $u(x_1)\in K[x_1]$ such that
\begin{equation}\label{eq4.1}
x_1^{p-1} = \delta(u(x_1))=u(x_1)-u(x_1+c)
\end{equation}
We have the following equations
\begin{equation}
  \left\{ \begin{aligned}
  \nonumber
 x_1^{p-1} = u(x_1)-u(x_1+c)~~~~~~~~~~~~~~~~~~~~~~~~~~~~~~~~\\
(x_1+c)^{p-1} = u(x_1+c)-u(x_1+2c)~~~~~~~~~~~~~~~~~~\\
\vdots~~~~~~~~~~~~~~~~~~~~~~~~~~~~~~~~~~~~~~~~~\\
(x_1+(p-1)c)^{p-1} = u(x_1+(p-1)c)-u(x_1+pc)
                          \end{aligned} \right.
  \end{equation}
by substituting $x_1$ with $x_1+jc$ for $0\leq j\leq p-1$ in equation \eqref{eq4.1}. Then we have
\begin{equation}\label{eq4.2}
\sum_{j=0}^{p-1}(x_1+jc)^{p-1} = 0
\end{equation}
by adding the above $p$ equations. By Fermat's little Theorem, we have $j^{p-1}=1 \operatorname{mod} p$ for all $1\leq j\leq p-1$. Let $x_1=0$ in equation \eqref{eq4.2}. Then we have that $(p-1)c^{p-1}=0$, which is a contradiction. Hence we have $x_1^{p-1}\notin \delta(I_0)$. Whence, $\delta(I_0)$ is not a Mathieu-Zhao space of $K[x_1]$.

If $c=0$, then $\phi(x_1)=x_1$, which is case $(1.3)$.

$(1.3)$ If $\phi(x_1)=\tilde{q}x_1$ for some $\tilde{q}\in K$, then $\delta(I_0)=K\{(1-\tilde{q})x_1, (1-\tilde{q}^2)x_1^2,\ldots,(1-\tilde{q}^k)x_1^k,\ldots\}$.

If $\tilde{q}$ is not a root of unity, then $\delta(I_0)$ is an ideal generated by $x_1$. Thus, $\delta(I_0)$ is a Mathieu-Zhao space of $K[x_1]$.

If $\tilde{q}=1$, then $\delta(I_0)=0$. Clearly, $\delta(I_0)$ is a Mathieu-Zhao space of $K[x_1]$.

If $\tilde{q}$ is a root of unity and $\tilde{q}\neq 1$, then there exists an integer $r\geq 2$ such that $\tilde{q}^r=1$. Let $S=\{n_j\in \mathbb{N}|~r\nmid n_j, j\geq 1\}$. Then  $\delta(I_0)$ is a $K$-subspace generated by $x_1^{n_j}$ for all $n_j\in S$. We claim that $\mathfrak{r}(\delta(I_0))=0$. Suppose that $\mathfrak{r}(\delta(I_0))\neq 0$. Then there exists nonzero polynomial $h(x_1)$ such that $h(x_1)\in \mathfrak{r}(\delta(I_0))$. Thus, $h^m(x_1)\in \delta(I_0)$ for all $m>>0$. Since $1\notin \delta(I_0)$, we have that $\operatorname{deg}h(x_1)=d_0\geq 1$. Since $\delta(I_0)$ is homogeneous, we have that $x_1^{d_0m}\in \delta(I_0)$ for all $m>>0$. That is, $d_0m\in S$ for all $m>>0$, which is a contradiction. Hence we have that $\mathfrak{r}(\delta(I_0))=0$. Therefore, $\delta(I_0)$ is a Mathieu-Zhao space of $K[x_1]$.\\

$(2)$ $(2.1)$ If $\phi(x_1)=c\neq 0$, then $\delta(I_0)=K\{x_1^i-c^i, (x_1-c)x_1^i, (x_1-c)x_1^{i+1},\ldots,\allowbreak(x_1-c)x_1^{i+k},\ldots\}$. Thus, we have that
 $$(x_1^i-c^i)^m = (x_1^i-c^i)(x_1^i-c^i)^{m-1}$$
 $$~~~~~~~~~~~~~~~~~~~~~~~~~~~~~~~~~~~~~~~~~=(x_1^i-c^i)(\sum_{l=0}^{m-1}\frac{(m-1)!}{(m-1-l)!l!}(x_1^i)^{m-1-l}(c^i)^l).$$                          It's easy to see that $(x_1^i-c^i)(x_1^i)^{m-1-l}(c^i)^l\in \delta(I_0)$ for all $0\leq l\leq m-1$ and for all $m\in \mathbb{N}^*$. Hence we have $(x_1^i-c^i)^m\in \delta(I_0)$ for all $m\geq 1$. We claim that $(x_1^i-c^i)^mx_1\notin \delta(I_0)$ for all $m\geq 2$. Suppose that $(x_1^i-c^i)^mx_1\in \delta(I_0)$ for some $m\geq 2$. Then we have that $(x_1^i-c^i)x_1\in \delta(I_0)$ because
$$(x_1^i-c^i)^mx_1=(x_1^i-c^i)(\sum_{l=0}^{m-1}\frac{(m-1)!}{(m-1-l)!l!}(x_1^i)^{m-1-l}(c^i)^l)x_1$$
and $(x_1^i-c^i)(x_1^i)^{m-1-l}(c^i)^l)x_1\in \delta(I_0)$ for all $0\leq l\leq m-2$ and for all $m\geq 2$. Since
$$c^i(x_1-c)=c(x_1^i-c^i)-[(x_1^i-c^i)x_1-(x_1-c)x_1^i],$$
we have that $x_1-c\in \delta(I_0)$, which is a contradiction. Thus, we have that \allowbreak $(x_1^i-c^i)^mx_1\notin \delta(I_0)$ for all $m\geq 2$. Hence $\delta(I_0)$ is not a Mathieu-Zhao space of $K[x_1]$.

If $c=0$, then $\delta$ is the identity map and $\delta(I_0)=I_0$. Clearly, $\delta(I_0)$ is a Mathieu-Zhao space of $K[x_1]$.

$(2.2)$ If $\phi(x_1)=x_1+c$ and $c\neq 0$, then $\delta(I_0)=K\{x_1^i-(x_1+c)^i, x_1^{i+1}-(x_1+c)^{i+1},\ldots\}$. If $i\leq p$, then $x_1^p-(x_1+c)^p=-c^p$. Thus, we have $1\in \delta(I_0)$. If $i>p$, then there exists $t\in \mathbb{N}^*$ such that $p^t\leq i\leq p^{t+1}$ and $x_1^{p^{t+1}}-(x_1+c)^{p^{t+1}}=-c^{p^{t+1}}$. Thus, we have $1\in \delta(I_0)$ in the situation. Then we have $x_1^{p-1}\notin \delta(I_0)$ by following the arguments of case $(1.2)$. Therefore, $\delta(I_0)$ is not a Mathieu-Zhao space of $K[x_1]$.

If $c=0$, then $\phi(x_1)=x_1$, which is case $(2.3)$.

$(2.3)$ If $\phi(x_1)=\tilde{q}x_1$ for some $\tilde{q}\in K$, then $\delta(I_0)=K\{(1-\tilde{q}^i)x_1^i, (1-\tilde{q}^{i+1})x_1^{i+1},\ldots,(1-\tilde{q}^{i+k})x_1^{i+k},\ldots\}$.

If $\tilde{q}$ is not a root of unity, then $\delta(I_0)$ is an ideal generated by $x_1^i$. Clearly, $\delta(I_0)$ is a Mathieu-Zhao space of $K[x_1]$.

If $\tilde{q}=1$, then $\delta(I_0)=0$. Clearly, $\delta(I_0)$ is a Mathieu-Zhao space of $K[x_1]$.

If $\tilde{q}$ is a root of unity and $\tilde{q}\neq 1$, then there exists an integer $r\geq 2$ such that $\tilde{q}^r=1$. Let $S=\{n_j\in \mathbb{N}|~r\nmid n_j, n_j\geq i, j\geq 1 \}$. Then $\delta(I_0)$ is a $K$-subspace generated by $x_1^{n_j}$ for all $n_j\in S$. Thus, we have $\mathfrak{r}(\delta(I_0))=0$ by following the arguments of case $(1.3)$. Therefore, $\delta(I_0)$ is a Mathieu-Zhao space of $K[x_1]$.
\end{proof}

{\bf{Acknowledgement}}:  The second author is very grateful to professor Wenhua Zhao for personal communications about the Mathieu-Zhao spaces. She is also grateful to the Department of Mathematics of Illinois State University, where this paper was partially finished, for hospitality during her stay as a visiting scholar.

\end{document}